\documentclass[aos,MSNbibl,seceqn,nameyear,dvips]{arximspdf}
\usepackage{mathrsfs}

\doi{10.1214/11-AOS892}
\volume{39}
\issue{4}
\pubyear{2011}
\firstpage{2007}
\lastpage{2020}

\makeatletter

\newcommand{\bfx}{\mathbf{x}}
\newcommand{\bfy}{\mathbf{y}}
\newcommand{\N}{{\mathbb N}}
\newcommand{\Nat}{{\mathbb N}}
\newcommand{\R}{{\mathbb R}}
\newcommand{\mcS}{\mathcal{S}}
\newcommand{\Po}{\operatorname{Po}}
\newcommand{\PP}{\operatorname{PP}}
\newcommand{\Obs}{\mathcal{S}_1}

\newcommand{\pr}{\operatorname{pr}}
\newcommand{\given}{\mid}

\newproclaim{example}{Example}[section]
\newproclaim{Method}{Method}

\makeatother

\begin{document}
\begin{frontmatter}

\title{On Bayes' theorem for improper mixtures\thanksref{T1}}
\runtitle{Bayes' theorem for improper mixtures}

\thankstext{T1}{Supported by NSF Grant DMS-09-06592.}

\begin{aug}
\author[A]{\fnms{Peter} \snm{McCullagh}\corref{}\ead[label=e1]{pmcc@galton.uchicago.edu}}
and
\author[A]{\fnms{Han} \snm{Han}\ead[label=e2]{han@galton.uchicago.edu}}
\runauthor{P. McCullagh and H. Han}
\affiliation{University of Chicago}
\address[A]{Department of Statistics\\
University of Chicago\\
5734 South University Ave\\
Chicago, Illinois 60637\\
USA\\
\printead{e1}\\
\hphantom{E-mail: }\printead*{e2}} 
\end{aug}

\received{\smonth{11} \syear{2010}}
\revised{\smonth{4} \syear{2011}}

%
\begin{abstract}
Although Bayes's theorem demands a prior that is a probability
distribution on the parameter space, the calculus associated with
Bayes's theorem sometimes generates sensible procedures from
improper priors, Pitman's estimator being a good example. However,
improper priors may also lead to Bayes procedures that are
paradoxical or otherwise unsatisfactory, prompting some authors to
insist that all priors be proper. This paper begins with the
observation that an improper measure on $\Theta$ satisfying
Kingman's countability condition is in fact a probability
distribution on the power set. We show how to extend a model in such
a way that the extended parameter space is the power set. Under an
additional finiteness condition, which is needed for the existence
of a sampling region, the conditions for Bayes's theorem are
satisfied by the extension. Lack of interference ensures that the
posterior distribution in the extended space is compatible with the
original parameter space. Provided that the key finiteness condition
is satisfied, this probabilistic analysis of the extended model may
be interpreted as a vindication of improper Bayes procedures derived
from the original model.
\end{abstract}

%

\end{frontmatter}

\section{Introduction}\label{secintro}
Consider a parametric model consisting of a family of probability
distributions $\{P_\theta\}$ indexed by the parameter
$\theta\in\Theta$. Each $P_\theta$ is a probability distribution on
the observation space $\Obs$, usually a product space such
as $\R^n$. In the parametric application of Bayes's theorem, the
family $\{P_\theta\}$ is replaced by a single probability
distribution $P_\pi(d\theta, dy) = P_\theta(dy) \pi(d\theta)$ on
the product space $\Theta\times\Obs$. The associated projections
are the prior $\pi$ on the parameter space
and the marginal distribution
\[
P_\pi(\Theta\times A) = \int_\Theta P_\theta(A) \pi(d\theta)
\]
for $A\subset\Obs$. To each observation, $y\in\Obs$ there corresponds
a conditional distribution $P_\pi(d\theta\given y)$, also called
the posterior distribution, on $\Theta$.

The joint distribution $P_\pi(d\theta, dy)$ has a dual
interpretation. The generative interpretation begins with $\theta$,
a random element drawn from $\Theta$ with probability
distribution $\pi$, the second component being distributed according
to the model distribution $Y\sim P_\theta$, now treated as a
conditional distribution given~$\theta$. In reverse order, the
inferential interpretation begins with the observational component
$Y\sim P_\pi(\Theta\times dy)$ drawn from the mixture distribution,
the parameter component being distributed as $\theta\sim
P_\pi(\cdot\given y)$ from the conditional distribution given
$Y=y$. The conditional distribution $P_\pi(\cdot\given y)$ tells us
how to select $\theta\in\Theta$ in order that the joint distribution
should coincide with the given joint distribution $P_\pi(d\theta, dy)$.

On the assumption that the marginal measure
$P_\nu(dy) = \int_\Theta P_\theta(dy) \nu(d\theta)$
is $\sigma$-finite,
formal application of the Bayes calculus with an improper prior~$\nu$
yields a posterior distribution $Q(d\theta\given y)$ satisfying
\[
P_\theta(dy) \nu(d\theta) = P_\nu(dy) Q(d\theta\given y)
\]
[\citet{EML}, \citet{ES}]. This factorization of the joint measure
yields a conditional law that is a probability distribution, in the
sense that $Q(\Theta\given y) = 1$. However, the joint measure is
not a probability distribution, so the factorization is not to be
confused with Bayes's theorem: it does not offer a probabilistic
interpretation of $Q(\cdot\given y)$ as a family of conditional
distributions generated by a joint probability distribution on the
product space. As a result, some authors reject the Kolmogorov axiom
of total probability, arguing instead for a nonunitary measure
theory for Bayesian applications [\citet{Hartigan},
\citet{Taraldsen}]. The goal of this paper is to show how an improper
prior may be accommodated within the standard unitary theory without
deviation from the Kolmogorov axioms. A probability space is
constructed from the improper measure in such a way that $Q(\cdot
\given y)$ admits a~probabilistic interpretation as a family of
conditional probability distributions given the observation.
Section~\ref{secBinary} shows that $\sigma$-finiteness is not
needed.

It would be inappropriate here to offer a review of the vast
literature on improper priors, most of which is not relevant to the
approach taken here. Nonetheless, a few remarks are in order. Some
statisticians clearly have qualms about the use of such priors,
partly because Bayes's theorem demands that priors be proper, partly
because the ``degree of belief'' interpretation is no longer
compelling, and partly because the formal manipulation of improper
priors may lead to inferential paradoxes of the sort discussed by
\citet{DZS}. \citet{L} argues correctly that strict adherence
to the
rules of probability requires all priors to be proper. Even though
the Bayes calculus often generates procedures yielding sensible
conclusions, he concludes that improper priors must be rejected.
Many statisticians, including some who interpret the prior as a
``degree of belief,'' are inclined to take a less dogmatic view. In
connection with Bernoulli trials, \citet{BS} (Section 5.2) comment
as follows. \textit{It is important to recognize, however, that this is
merely an approximation device and in no way justifies} [\textit{the improper
limit $\theta^{-1}(1-\theta)^{-1}$}] \textit{as having any special
significance as a representation of} ``\textit{prior ignorance}.'' In
subsequent discussion in Section 5.4, they take a more pragmatic
view of a reference prior as a mathematical tool generating a
reference analysis by the Bayes calculus.

The purpose of this note is to offer a purely probabilistic
interpretation of an improper prior, in agreement with Lindley's
thesis but not with his conclusion. The interpretation that removes
the chief mathematical obstacle
is that an improper measure on $\Theta$ is a probability
distribution on the set of subsets of $\Theta$. A proper prior
determines a random \textit{element} $\theta\in\Theta$ with
distribution $\pi$, whereas an improper prior $\nu$ determines a
random \textit{subset}, a~countable collection
$\{\theta_i\}$ distributed as a Poisson process with mean
measure~$\nu$. In the product space $\Theta\times\Obs$, the proper joint
distribution $P_\pi$ determines a random element $(\theta, Y)$, whereas
the improper distribution $P_\nu$ determines a random subset $Z\subset
\Theta\times\Obs$, a countable collection of ordered pairs $Z =
\{(\theta_i, Y_i)\}$. An observation on a point process consists of
a sampling region $A\subset\Obs$ together with the set $\bfy= Y\cap
A$ of events that occur in~$A$. It is critical that the sampling
region be specified in such a way that $Y\cap A$ is finite, a
condition that puts definite limits on $\nu$ and on the set of
sampling schemes. Having done so, we obtain the conditional
distribution given the observation. The standard Bayesian argument
associates with each point $y\in\Obs$ a probability distribution
on~$\Theta$: the point process argument associates with each finite
subset $\bfy\subset A$ a probability distribution on
$\Theta^{\#\bfy}$. Despite this fundamental distinction, certain
aspects of the conditional distribution are in accord with the
formal application of the Bayes calculus, treating the mixture as if
it were a model for a random element rather than a random subset.

\section{Conditional distributions}\label{secconditional}
Consider a Poisson process with mean measure $\mu$ in the product
space $\mcS= \mcS_0\times\mcS_1$. Existence of the process is guaranteed
if the singletons of $\mcS$ are contained in the $\sigma$-field, and
$\mu$ is a countable sum of finite measures, that is,
%
%
\begin{equation}\label{kfc}
\mu= \sum_{n=1}^\infty\mu_n \qquad\mbox{where } \mu_n(\mcS) <
\infty.
\end{equation}
Kingman's countability condition, also called weak finiteness
[Kingman\break (\citeyear{King})], is the natural condition for existence because it
implies that the marginal measures $\mu_0(B) = \mu(B\times\mcS_1)$ for
$B\subset\mcS_0$ and $\mu_1(A) = \mu(\mcS_0\times A)$ for $A\subset\mcS_1$
are countable. Consequently, the projected processes exist and are
also Poisson.

Unlike $\sigma$-finiteness, countability does not imply the
existence of a subset $A\subset\mcS$ such that $0 < \mu(A) < \infty$.
If such a set exists, the process is said to be \textit{observable
on~$A$}. For example, the measure taking the value $\infty$ on
subsets of positive Lebesgue measure in $\R$ and zero otherwise is
countable, but the process is not observable on any subset.
Sigma-finiteness is a stronger condition, sufficient for existence
but not necessary, and not inherited by the projected marginal
measures [\citet{King}].

The symbol $Z\sim\PP(\mu)$ denotes a Poisson point
process, which is a random subset $Z\subset\mcS$ such that for each
finite collection of disjoint subsets $A_1,\ldots, A_n$ of~$\mcS$, the
random variables $\#(Z\cap A_1), \ldots, \#(Z\cap A_n)$ are
distributed independently according to the Poisson distribution
$\#(Z\cap A_j) \sim\Po(\mu(A_j))$. In much of what follows, it is
assumed that $\mu(\mcS) = \infty$, which implies that $\#Z
\sim\Po(\infty)$ is infinite with probability one, but countable on
account of~(\ref{kfc}). Since $Z$ is countable and $\mcS$ is a product
set, we may label the events
\[
Z = (X, Y) = \{(X_i, Y_i) \dvtx i=1,2,\ldots\},
\]
where $X\subset\mcS_0$ is a Poisson process with mean measure $\mu_0$
and $Y\subset\mcS_1$ is a Poisson process with mean measure $\mu_1$.
The notation $Z = (X, Y)$ implies that $X\subset\mcS_0$ and
$Y\subset\mcS_1$ are countable subsets whose elements are in a~specific 1--1 correspondence.

To say what is meant by an \textit{observation} on a point process, we
must first establish the sampling protocol, which is a test set or
sampling region $A\subset\mcS_1$ such that $\mu_1(A) < \infty$.
In this scheme, $\mcS_0$ is the domain of inference, so $X$ is not
observed. The actual observation is the test set $A$ together with
the random subset $\bfy= Y\cap A$, which is finite with probability
one. Although we refer to $\mcS_1$ as the ``space of observations,'' it
must be emphasized that an observation is not a random
\textit{element} in~$\mcS_1$, but a finite random \textit{subset}
$\bfy\subset A\subset\mcS_1$, which could be empty.

The distinction between a point process and an observation on the
process is the same as the distinction between an infinite process
and an observation on that process. An infinite process is a
sequence of random variables $Y=(Y_1, Y_2,\ldots)$ indexed by the
natural numbers, that is, a random function $Y\colon\Nat\to\R$. An
observation consists of a sample, a \textit{finite} subset $A \subset
\Nat$, together with the response values $Y[A]$ for the sampled
units. Likewise, a point process is a random subset considered as a
random function $Y\colon\mcS_1 \to\{0,1\}$ indexed by the
domain $\mcS_1$. An observation consists of a sample or sampling
region $A\subset\mcS_1$ together with the restriction $Y[A] = Y\cap A$
of the process to the sample. Usually $A$ is not finite or even
countable, but the observation is necessarily finite in the sense
that $\#(Y\cap A) < \infty$.

Whether we are talking of sequences or point processes, the domain
of inference is not necessarily to be interpreted as a parameter
space: in certain applications discussed below, the observation
space consists of finite sequences in $\mcS_1=\R^n$, and\vadjust{\goodbreak}
$\mcS_0 = \R^\infty$ is the set of subsequent trajectories. In this sense,
predictive sample-space inferences are an integral part of the
general theory (Section \ref{secnonpara}).

We focus here on inferences for the $X$-values associated with the
events $\bfy=Y\cap A$ that occur in the sampling region, that is, the
subset
\[
\bfx= X[A] = \{X_i \dvtx Y_i \in A\} = \{X_i \dvtx Y_i\in\bfy\}
\]
in 1--1 correspondence with the observation $\bfy$. In this formal
sense, an inference is a rule associating with each finite subset
$\bfy\subset A$ a probability distribution on $\mcS_0^{\#\bfy}$.

Clearly, if $\bfy$ is empty, $\bfx$ is also empty, so the
conditional distribution is trivial, putting probability one on the
event that $\bfx$ is empty. Without loss of generality, therefore,
we assume that $0 < \mu_1(A) < \infty$, that $m=\#\bfy$ is positive
and finite, and that the events are labeled $(Y_1,\ldots, Y_m)$ by a~uniform
random permutation independent of $Z$. Given $\#\bfy=m$, the
pairs $(X_1,Y_1),\ldots, (X_m, Y_m)$ are independent and identically
distributed random variables with probability density $\mu(dx
dy)/\mu_1(A)$ in $\mcS_0\times A$. Thus, the conditional joint density
given $Y\cap A = \bfy$ is equal to
%
%
\begin{equation}\label{conditionaldistribution}
p(d\bfx\given\bfy) = \prod_{i=1}^m \frac{\mu(dx_i
dy_i)}{\mu_1(dy_i)}
=\prod_{i=1}^m \mu(dx_i \given y_i),
\end{equation}
where $\mu(dx \given y)$ is the limiting ratio $\mu(dx\times
dy)/\mu_1(dy)$ as $dy\downarrow\{y\}$.

The key properties of this conditional distribution are twofold,
\textit{conditional independence} and \textit{lack of interference}.
First, the random variables $X_1,\ldots, X_m$ are conditionally
independent given $Y\cap A=\bfy$. Second, the conditional
distribution of $X_i$ given $\bfy$ depends only on $Y_i$, not on the
number or position of other events in $A$. For example, if two or
more events occur at the same point ($Y_i = Y_j$) the random
variables $X_i, X_j$ are conditionally independent and identically
distributed given $\bfy$. The test set determines the events on
which predictions are made, but beyond that it has no effect. In
particular, if $m=1$, the conditional density of $X$ is $p(dx \given
y) \propto\mu(dx \given y)$ regardless of the test set.

The observability assumption $\mu_1(A) < \infty$ is not made out of
concern for what might reasonably be expected of an observer in the
field. On the contrary, finiteness is essential to the mathematical
argument leading to (\ref{conditionaldistribution}). If the number
of events were infinite, countability implies that the values can be
labeled sequentially $y_1, y_2,\ldots$ in 1--1 correspondence with
the integers. Countability does not imply that they can be labeled
in such a way that the infinite sequence is exchangeable.
As a result, the factorization
(\ref{conditionaldistribution})
fails if $\#\bfy=\infty$.

The remark made above, that the test set has no effect on
inferences, is correct but possibly misleading. Suppose that $0 < m
< \infty$ and that the observation consists of that information
alone without recording the particular values. If $\mu_1(A) = 0$ or
$\mu_1(A) = \infty$, no inference is possible beyond\vadjust{\goodbreak} the fact that
the model is totally incompatible with the observation. If the
marginal measure is finite on $A$, the conditional density is such
that the components of $X[A]$ are independent and identically
distributed with density $\mu(dx\times A) / \mu_1(A)$, which does
depend on the choice of test set. In the context of parametric
mixture models with $\Theta\equiv\mcS_0$, each sequence with
distribution $P_\theta$ has probability $P_\theta(A)$ of being
recorded. Thus, before observation begins, the restriction to
$A\subset\mcS_1$ effectively changes the measure to $P_\theta(A)
\nu(d\theta)$, which is finite on $\Theta$, but depends on the
choice of $A$.

\vspace*{3pt}\section{Improper mixtures}\label{secImproper}
Consider a parametric statistical model consisting of a family of
probability distributions $\{P_\theta\colon\theta\in\Theta\}$ on
the observation space $\mcS_1$, one distribution for each point
$\theta$ in the parameter space. Each model distribution determines
a random \textit{element} $Y \sim P_\theta$. A probability
distribution $\pi$ on $\Theta$ completes the Bayesian specification,
and each Bayesian model also determines a random element
$(\theta,Y)\in\Theta\times\mcS_1$ distributed as $\pi(d\theta)
P_\theta(dy)$. The observational component is a random element $Y\in
\mcS_1$ distributed as the mixture $Y\sim P_\pi$, and the conditional
distribution given $Y=y$ is formally the limit of $\pi(d\theta)
P_\theta(dy)/P_\pi(\Theta,dy)$ as $dy\downarrow\{y\}$.

A countable measure $\nu$ such that $\nu(\Theta) = \infty$ does not
determine a \textit{random element} $\theta\in\Theta$, but it does
determine an infinite \textit{random subset} $X\subset\Theta$.
Furthermore, the joint measure $\nu(d\theta) P_\theta(dy)$ is
countable, so there exists a~random subset $Z=(X,Y)
\subset\Theta\times\mcS_1$, distributed according to the Poisson
process with mean measure $\nu(d\theta) P_\theta(dy)$. If this
interpretation is granted, it is necessary first to specify the
sampling region $A\subset\mcS_1$, in such a way that $P_\nu(A) <
\infty$ to ensure that only finitely many events $\bfy=Y\cap A$
occur in $A$. To each observed event $Y_i\in\bfy$, there corresponds
a parameter point $\theta_i\in X[A]$ such that $(\theta_i, Y_i) \in
Z$. Parametric inference consists in finding the joint conditional
distribution given $Y\cap A = \bfy$ of the particular subset of
parameter values $\theta_1,\ldots, \theta_m$ corresponding to the
events observed.

This probabilistic interpretation forces us
to think of the parameter and the observation in a
collective manner, as sets rather than points. Taken
literally, the improper mixture is not a model for a random element
in~\mbox{$\Theta\times\mcS_1$}, but a model for a random subset
$Z= (X,Y)\subset\Theta\times\mcS_1$. If $\nu(\Theta) < \infty$,
as in a proper mixture,
it is sufficient to take $A=\mcS_1$ and to record the entire
subset $\bfy\subset\mcS_1$, which is necessarily finite.
However, if $\nu(\Theta) = \infty$, it is necessary to sample the process
by first establishing a test set $A\subset\mcS_1$ such that
$P_\nu(A) < \infty$, and then listing the finite set of values
$\bfy=Y\cap A$ that occur in $A$. Generally speaking, this
finiteness condition rules out many sampling schemes that might
otherwise seem reasonable. In the special case where $\#\bfy=1$,
$X[A]$ is a random subset consisting of a single point, whose
conditional density at $x\in\Theta$ is
%
%
\begin{equation}
\pr(X[A] \in dx \given\bfy=\{y\}) = \frac{\nu(dx)
p_x(y)}{\int_\Theta p_\theta(y) \nu(d\theta)},
\end{equation}
where $p_\theta(y)$ is the density of $P_\theta$ at $y$. The
finiteness condition on $A$ ensures that the integral in the
denominator is finite, and the occurrence of an event at~$y$ implies
that $P_\nu$ assigns positive mass to each open neighborhood of $y$.

Provided that $0 < P_\nu(A) < \infty$, this purely probabilistic conclusion
may be interpreted as a vindication of the formal Bayes calculation
associated with an improper prior. However, the two versions of
Bayes's theorem are quite different in logical structure; one
implies a single random element, the other infinitely many.
Accordingly, if a statistical procedure is to be judged by a
criterion such as a conventional loss function, which presupposes a
single observation and a single parameter, we should not expect
optimal results from a probabilistic theory that demands multiple
observations and multiple parameters. Conversely, if the procedure
is to be judged by a criterion that allows for multiple sequences
each with its own parameter, we should not expect useful results
from a probabilistic theory that recognizes only one sequence and
one parameter. Thus, the existence of a joint probability model
associated with an improper prior does not imply optimality in the
form of coherence, consistency or admissibility. For example, in the
MANOVA example of
\citet{ES}, the Poisson point process interpretation
yields the classical posterior, which is incoherent in de Finetti's
sense and is strongly inconsistent in Stone's sense.

The observability condition implies that the restriction of $P_\nu$
to $A$ is finite, and hence trivially $\sigma$-finite. The role of
the finiteness condition is illustrated by two examples in
Sections \ref{secGaussian} and \ref{secBinary}. For the Gaussian
model, $P_\nu$ is countable for every $n\ge0$ and $\sigma$-finite
for $n\ge2$, which guarantees the existence of a~sampling region if
$n\ge2$. For the Bernoulli model, $P_\nu$ is countable for each
$n\ge0$ but not $\sigma$-finite for any $n$. Nonetheless, the
finiteness condition for observability is satisfied by certain
subsets $A\subset\{0,1\}^n$ for $n\ge2$.

\section{Gaussian point process}\label{secGaussian}
\subsection{Parametric version}
Consider the standard model for a Gaussian sequence with independent
$N(\theta, \sigma^2)$ components. Let $p$ be a given real number,
and let the prior measure be $\nu(d\theta \,d\sigma) = d\theta
\,d\sigma/\sigma^p$ on the parameter space $\Theta= \R\times\R^+$.
For all $p$, both $\nu$ and the joint measure on $\Theta\times\R^n$
satisfy the countability condition. Consequently a Poisson point
process $Z=(X,Y)\subset\Theta\times\R^n$ exists in the product
space. For $n > 2-p$, the marginal measure $P_\nu$ has a density in
$\R^n$
%
%
\begin{equation}\label{intensity}
\lambda_n(y) = \frac{\Gamma((n+p-2)/2) 2^{(p-3)/2} \pi^{-(n-1)/2}
n^{-1/2}}
{(\sum_{i=1}^n (y_i - \bar y)^2 )^{(n+p-2)/2}},
\end{equation}
which is finite at all points $y\in\R^n$ except for the diagonal
set. Provided that $n\ge2$ and $n > 2-p$, there exists in $\R^n$ a
subset $A$ such that \mbox{$P_\nu(A) < \infty$}, which serves as the region
of observation. In fact, these conditions are sufficient for
$\sigma$-finiteness in this example. To each observation $\bfy=Y\cap
A$ and to each event $y\in\bfy$, there corresponds a conditional
distribution on $\Theta$ with density
\[
p(\theta, \sigma\given Y\cap A=\bfy, y\in\bfy) =
\phi_n(y; \theta, \sigma) \sigma^{-p} / \lambda_n(y),
\]
where $\phi_n(y; \theta, \sigma)$ is the Gaussian density at $y$ in
$\R^n$. The conditional distribution (\ref{conditionaldistribution})
of the parameter subset $X[A]\subset\Theta$ given $Y\cap A=\bfy$ is
a product of factors of this type, one for each of the events
in $\bfy$. It should be emphasized here that the information in the
conditioning event is not simply that $\bfy\subset Y$, but also that
$Y$ contains no other events in $A$.

\subsection{Nonparametric version}\label{secnonpara}
Let $\mathbb{N}$ be the set of natural numbers, and let
$\mcS=\R^{\mathbb N}$ be the collection of real-valued sequences,
\[
\mcS=\R^{\mathbb{N}}=\{y=(y_1,y_2,\ldots) \dvtx y_i\in\mathbb{R}, i
\in
\mathbb{N}\}
\]
with product $\sigma$-field $\mathscr{R}^{\N}$. We construct
directly in this space a Poisson process $Z\subset\mcS$ whose mean
measure $\Lambda$ is uniquely determined by its finite-dimensional
projections $\Lambda_n$ with density (\ref{intensity}). By their
construction, these measures are finitely exchangeable and satisfy
the Kolmogorov consistency condition $\Lambda_{n+1}(A\times\R) =
\Lambda_n(A)$ for each integer $n\ge0$ and $A\in\mathscr R^n$. In
keeping with the terminology for random sequences, we say that the
point process $Z\sim\PP(\Lambda)$ is infinitely exchangeable if each
$\Lambda_n$ is finitely exchangeable.\looseness=1

Let $\mcS=\mcS_1\times\mcS_0$, where $\mcS_1=\R^n$ is the projection onto
the first $n$ coordinates, and $\mcS_0\cong\mcS$ is the complementary
projection onto the subsequent coordinates. Each event $z\in Z$ is
an ordered pair, so we write $Z=(Y, X)\subset\mcS$ as a countable set
of ordered pairs $(Y_i, X_i)$ in which the marginal process
$Y\subset\mcS_1$ is Poisson with parameter $\Lambda_n$, and
$X\sim\PP(\Lambda)$ has the same distribution as $Z$. Provided that
the set $A\subset\mcS_1$ has finite $\Lambda_n$-measure, the
observation $\bfy= Y\cap A$ is finite. To each event $y\in\bfy$,
there corresponds an event $z=(y,x)\in Z$, so that $y=(z_1,\ldots,
z_n)$ is the initial sequence, and $x=(z_{n+1},\ldots)$ is the
subsequent trajectory. The conditional distribution~(\ref{conditionaldistribution}) is such that the subsequent
trajectories $X[A]$ are conditionally independent and
noninterfering given $Y\cap A = \bfy$. For each event
$y\in\bfy$, the $k$-dimensional joint density at $x = (x_1, \ldots,
x_k)$ of the subsequent trajectory is
%
%
\begin{equation}\label{prediction}
p(dx \given Y\cap A=\bfy, y\in\bfy) =
\frac{\lambda_{n+k} (y, x) \,dx} {\lambda_n(y)},
\end{equation}
which is the $k$-dimensional exchangeable Student $t$ density
[\citet{KN}, page 1] on $\nu=n+p-2>0$ degrees of freedom.

For any\vspace*{1pt} continuous location-scale model with finite $p$th moment
and improper prior density proportional to $d\mu \,d\sigma/\sigma^p$
with $p>0$, the initial segment $Y\subset\R^2$ is a Poisson process
with intensity
\[
\lambda_2(y) \propto\frac{1} {|y_1 - y_2|^p}.
\]
Otherwise if $p\leq0$ the initial segment of length $n>2-p$ is a
Poisson process with intensity
\[
\lambda_n(y) \propto\frac{1}
{(\sum_{i=1}^n (y_i - \bar y)^2 )^{(n+p-2)/2}}.
\]
%
The prescription (\ref{prediction})
extends each event $y\in Y$ to an infinite random sequence in such a
way that the set of extended sequences $Z\subset\R^{\mathbb{N}}$ is
a Poisson process with mean measure $\Lambda$. Given $Y\subset\R^2$,
these extensions are conditionally independent, noninterfering, and
each extension is an exchangeable sequence. In the Gaussian case,
(\ref{prediction}) is equivalent to the statement that each initial
sequence with $s_n\neq0$ is extended according to the recursive
Gosset rule 
\[
y_{n+1} = \bar y_{n} + s_{n} \varepsilon_{n}
\sqrt{\frac{n^2-1}{n(n+p-2)}},
\]
where $\bar y_n, s_n^2$ are the sample mean and variance of the
first $n$ components, and $\varepsilon_{n} \sim t_{n+p-2}$ has
independent components. The resulting extension is an exchangeable
sequence whose $k$-dimensional joint density at $(y_3, \ldots,
y_{k+2})$ is $\lambda_{k+2}(y_1,\ldots, y_{k+2}) / \lambda_2(y_1,
y_2)$.

The Gosset extension is such that the sequence $(\bar y_n, s_n^2)$
is Markovian and has a limit. Given a single sequence in the
sampling region, the joint distribution of the limiting random
variables $(\bar y_\infty, s_\infty)$ is
\[
p(\bar y_\infty, s_\infty\given Y\cap A=\bfy, y\in\bfy) =
\phi_n(y; \bar y_\infty,
s_\infty) s_{\infty}^{-p} / \lambda_n(y),
\]
which is the posterior density on $\Theta$ as computed by the Bayes
calculus with improper prior.

\vspace*{3pt}\section{Cauchy sequences}\label{secCauchy}
Consider the standard model for a Cauchy sequence having independent
components with parameter $\theta\in\R\times\R^+$. For $p > 0$, the
prior measure $\nu(d\theta) = d \theta_1 \, d\theta_2/\theta_2^p$
satisfies the countability condition, which implies that a Poisson
process $X = (X, Y)\subset\Theta\times\R^n$ with mean measure
$P_\nu$ exists in the product space. If $0<p < n$ and $n\geq2$, the
marginal measure in $\R^n$ has a density which is finite at all
points $y\in\R^n$ whose components are distinct. The density
satisfies the recurrence formula
\[
\lim_{y_n\to\pm\infty} \pi y_n^2 \lambda_{n,p}(y_1,\ldots, y_n) =
\lambda_{n-1, p-1}(y_1,\ldots, y_{n-1}).
\]
For integer $p\geq2$, the density is
%
%
\begin{equation}
\openup1\jot\lambda_{n,p}(y) = \cases{
\displaystyle {\frac{(-1)^{(n-p+1)/2}}{\pi^{n-2} 2^{n-p+1}}}
\sum_{r\neq s} { \frac{|y_s-y_r|^{n-p}}{d_r d_s}}, \qquad (n-p)
\mbox{ odd};\cr
\displaystyle {\frac{(-1)^{(n-p)/2}}{\pi^{n-1} 2^{n-p}}}
\sum_{r\neq s} { \frac{(y_s-y_r)^{n-p} \log|y_s-y_r|}{d_r d_s}}, \cr
\hspace*{169pt}(n-p)
\mbox{ even;}}
\end{equation}
where $d_r=\prod_{t\neq r}(y_t-y_r)$. For example, $\lambda_{2,1}(y)
= {1/(2|y_1-y_2|)}$
and
\[
\lambda_{3,2}(y) =
\frac{1}{2\pi|(y_1-y_2)(y_2-y_3)(y_1-y_3)|}.
\]
Spiegelhalter (\citeyear{Spi85}) established the same formula for $p=1$
in equation (2.2).

For $n> p$, there exists a subset $A\subset\R^n$ such that
$\Lambda_n(A) < \infty$, which serves as the region of observation.
The Poisson process determines a probability distribution on finite
subsets $\mathbf{y} \subset$ $A$, and to each point $y\in\mathbf{y}$ it also
associates a conditional distribution on $\Theta$ with density
%
%
\begin{equation}\label{C-posterior}
\frac{P_\nu(d\theta\times dy)} {\Lambda_n(dy)} = \frac{f_n(y;
\theta) \theta_2^{-p}} {\lambda_n(y)},
\end{equation}
where $f_n(y; \theta)$ is the Cauchy density at $y\in\R^n$.

In the nonparametric version with $\Theta$ replaced by $\R^k$, the
conditional distribution extends each point $y\in A$ to a sequence
$(y, X)\in\R^{n+k}$, with conditional density $X\sim
\lambda_{n+k}(y, x)/\lambda_n(y)$. The extension is infinitely
exchangeable. The tail trajectory of the infinite sequence is such
that, if $T_{k}\colon\R^k\to\Theta$ is Cauchy-consistent,
$T_{n+k}(y, X)$ has a limit whose density at $\theta\in\Theta$ is
(\ref{C-posterior}).

\vspace*{-3pt}\section{Binary sequences}\label{secBinary}
Consider the standard model for a Bernoulli sequence with parameter
space $\Theta=(0,1)$. The prior measure $\nu(d\theta) =
d\theta/(\theta(1-\theta))$ determines a Poisson process with
intensity $\theta^{n_1(y)-1} (1-\theta)^{n_0(y)-1}$ at $(y,\theta)$
in the product space $\mcS_1\times\Theta$. Here $\mcS_1 = \{0,1\}^n$ is
the space of sequences of length $n$, $n_0(y)$ is the number of
zeros and $n_1(y)$ is the number of ones in $y$. The marginal
measure on the observation space is
\[
\Lambda_n(\{y\}) = \cases{
\Gamma(n_0(y)) \Gamma(n_1(y)) / \Gamma(n), & \quad$n_0(y),
n_1(y) > 0$,\cr
\infty, & \quad otherwise,}
\]
which is countable but not $\sigma$-finite. Any subset
$A\subset\mcS_1$ that excludes the zero sequence and the unit sequence
has finite measure and can serve as the region of observation. Given
such a set and the observation $\bfy=Y\cap A$ recorded with
multiplicities, the conditional distribution
(\ref{conditionaldistribution}) associates with each $y\in\bfy$ the
beta distribution
\[
P_\nu(\theta\given Y\cap A=\bfy, y\in\bfy) =
\frac{\theta^{n_1(y) - 1} (1-\theta)^{n_0(y)-1} \Gamma(n)}
{\Gamma(n_1(y)) \Gamma(n_0(y))}
\]
on the parameter space.\vadjust{\goodbreak}

As in the preceding section, we may bypass the parameter space and
proceed directly by constructing a Poisson process with mean measure
$\Lambda$ in the space of infinite binary sequences. The values
assigned by $\Lambda$ to the infinite zero sequence and the infinite
unit sequence are not determined by~$\{\Lambda_n\}$, and can be set
to any arbitrary value, finite or infinite. Regardless of this
choice, (\ref{conditionaldistribution}) may be used to predict the
subsequent trajectory of each of the points $\bfy=Y\cap A$ provided
that $\Lambda_n(A) < \infty$. In particular, the conditional
distribution of the next subsequent component is
\[
\pr(y_{n+1} = 1 \given Y\cap A=\bfy, y\in\bfy) = n_1(y)/n.
\]
This is the standard P\'{o}lya urn model [\citet{D}] for which the
infinite average of all subsequent components is a beta random
variable with parameters $(n_0(y), n_1(y))$, in agreement with the
parametric analysis.

\vspace*{-3pt}\section{Interpretation}
The point-process interpretation of an improper measure on $\Theta$ forces
us to think of the parameter in a collective sense as a random
subset rather than a random point.
One interpretation is that a proper prior is designed for a specific
scientific problem whose goal is the estimation of a particular parameter
about which something may be known, or informed guesses can be made.
An improper mixture is designed for a generic class of problems,
not necessarily related to one another scientifically,
but all having the same mathematical structure.
Logistic regression models, which are used for many purposes
in a wide range of disciplines, are generic in this sense.
In the absence of a specific scientific context, nothing can be
known about the parameter, other than the fact that there
are many scientific problems of the same mathematical type,
each associated with a different parameter value.
In that wider sense of a generic mathematical class,
it is not unnatural to consider a broader framework
encompassing infinitely many scientific problems,
each with its own parameter.
The set of parameters is random but not indexed in an exchangeable way.

A generic model may be tailored to a specific scientific application
by coupling it with a proper prior distribution $\pi$ that is deemed
relevant to the scientific context. If there is broad agreement
about the model and the relevance of $\pi$ to the context,
subsequent calculations are uncontroversial. Difficulties arise when
no consensus can be reached about the prior. According to one
viewpoint, each individual has a personal prior or belief; Bayes's
theorem is then a recipe for the modification of personal beliefs
[\citet{BS}, Chapter 2]. Another line of argument calls for a panel
of so-called experts to reach a consensus before Bayes's theorem can
be used in a mutually agreeable fashion [\citet{Weerhandi},
\citet{Genest}]. A third option is to use proper but flat or
relatively uninformative priors. Each of these options demands a
proper prior on $\Theta$ in order that Bayes's theorem may be used.

This paper offers a fourth option by showing that
it is possible to apply Bayes's theorem to the generic model.
Rather than forcing the panel to reach\vadjust{\goodbreak} a proper consensus,
we may settle for an improper prior as a countable sum of
proper, and perhaps mutually contradictory, priors generated
by an infinite number of experts.
Although Bayes's theorem can be used,
the structure of the theorem for an improper mixture is not
the same as the structure for a proper prior.
For example, improper Bayes estimators need not be admissible.

Finiteness of the restriction of the measure to the sampling region is
needed in our argument. If the restriction to the sampling region is
$\sigma$-finite, we may partition the region into a countable family
of disjoint subsets of finite measure, and apply the extension subset
by subset. The existence of a Poisson point process on the sampling
region is assured by Kingman's superposition theorem. Lack of
interference implies that these extensions are
mutually consistent, so there is no problem dealing with
such $\sigma$-finite restrictions. This is probably not necessary from
a statistical perspective, but it does not create any mathematical
problems because the extension does not depend on the choice of the
partition of the region.

\vspace*{-3pt}\section{Marginalization paradoxes} \label{secMP}
$\!\!\!\!\!$The unBayesian characteristic of an impro\-per prior distribution is
highlighted by the marginalization paradoxes discus\-sed by \citet{SD}
and by \citet{DZS}.
In~the following example from \citet{SD}, the formal marginal
posterior distribution calculated by two methods demonstrates the
inconsistency.

\begin{example}\label{ex3}
The observation consists of two independent exponential random
variables
$X\sim\mathcal{E}(\theta\phi)$ and $Y\sim\mathcal{E}(\phi)$, where
$\theta$ and $\phi$ are unknown parameters. The parameter of
interest is the ratio $\theta$.
\end{example}
\begin{Method}\label{method1} The joint density is
\[
\pr(dx,dy\given\theta,\phi)=\theta\phi^2 e^{-\phi(\theta x +y)}
\,dx\,
dy.
\]
Given the improper prior distribution $\pi(\theta) \,d\theta
d \phi$, the marginal posterior distribution for $\theta$ is
%
%
\begin{equation}\label{post3}
\pi(\theta\given x,y)\propto\frac{\pi(\theta) \theta}{(\theta
x+y)^3}\propto\frac{\pi(\theta) \theta}{(\theta+ z)^3},
\end{equation}
where $z=y/x$.
\end{Method}
\begin{Method}\label{method2}
Notice that the posterior distribution depends on
$(x,y)$ only through $z$. For a given $\theta$, $z/\theta$ has an
$F_{2,2}$ distribution, that is,
\[
\pr(z\given\theta)\propto\frac{\theta}{(\theta+z)^2}.
\]
Using the implied marginal prior $\pi(\theta) \,d\theta$, as if it
were the limit of a sequence of proper priors, we obtain
%
%
\begin{equation}\label{post4}
\pi(\theta\given z)\propto\frac{\pi(\theta) \theta}{(\theta+ z)^2},\vadjust{\eject}
\end{equation}
which differs from (\ref{post3}). It has been pointed out by Dempster
and in the author's rejoinder [\citet{DZS}], that no choice of $\pi
(\theta)$ could bring the two analyses into agreement.

From the present viewpoint, the improper prior determines a random
subset of the parameter space and a random subset of the observation
space~$(\R^+)^2$. Under suitable conditions on $\pi$, the bivariate
intensity
\[
\lambda(x,y) 
= 2\int_0^\infty\frac{\theta \pi(\theta) \,d\theta}{(\theta x +
y)^3}
\]
is finite on the interior of the observation space, so the bivariate
process is observable. Equation (\ref{conditionaldistribution})
associates with each event $(x,y)$ the conditional distribution
(\ref{post3}) in agreement with the formal calculation by Method \ref{method1}.
Each event $(x, y)$ determines a ratio $z=y/x$, and the set of
ratios is a~Poisson point process in $(0,\infty)$. However, the
marginal measure is such that $\Lambda_z(A)=\infty$ for sets of
positive Lebesgue measure, and zero otherwise. This measure is
countable, but the marginal process is not observable. Thus,
conclusion~(\ref{post4}) deduced by Method \ref{method2} does not follow from
(\ref{conditionaldistribution}), and there is no contradiction.

Conversely, if the prior measure $\pi(d\theta)\rho(d\phi)$ is multiplicative
with $\rho(\R^+) < \infty$ and $\pi$ locally finite,
the marginal measure on the observation space is such that
\[
P_\nu(\{ a < x/y < b\}) < \infty
\]
for $0 < a < b < \infty$.
Thus, the ratio $z=x/y$ is observable,
and the conditions for Method \ref{method2} are satisfied.
The point process model associates with each ratio $0 < z < \infty$
the conditional distribution with density
\[
\pi(\theta, \phi\given z) \propto
\frac{\rho(\phi) \theta} {(\theta+z)^2}
\]
in agreement with (8.2).
However, the conditional distribution given $(x, y)$
\[
\pi(\theta, \phi\given(x,y)) \propto
\theta\phi^2 e^{-\phi(\theta x + y)} \pi(\theta) \rho(\phi)
\]
is such that the marginal distribution of $\theta$ given $(x,y)$
is not a function of $z$ alone.
Once again, there is no conflict with (8.2).

All of the other marginalization
paradoxes in \citet{DZS} follow the same pattern.

\citet{Jaynes} asserts that
``an improper pdf has meaning only as the limit of a well-defined
sequence of proper pdfs.''
On this point, there seems to be near-universal agreement,
even among authors who take diametrically opposed views
on other aspects of the marginalization paradox
[\citet{AK}, \citet{DZS} and \citet{Wall}].
No condition of this sort occurs in the point-process theory.
However, a sequence of measures $\mu_n$ such that $\mu_n(A) < \infty$,
each of which assigns a conditional distribution (2.2) to every $y\in A$,
may have a weak limit $\mu_n\to\mu$ such that $\mu(A) = \infty$
for which no conditional distribution exists.
\end{Method}


%

%
\printaddresses

\end{document}